\documentclass[
10pt, 
a4paper, 
oneside, 
]{article}

\usepackage{amssymb}
\usepackage{amsmath}
\usepackage{hyperref}

\usepackage[capitalise,noabbrev]{cleveref}
\usepackage{MnSymbol}
\usepackage{xcolor}

\usepackage{graphicx} 

\newtheorem{remark}{Remark}
\usepackage{algorithm}
\usepackage{algorithmic}

\providecommand{\keywords}[1]{\textbf{Key words. } #1}

\providecommand{\amssubj}[1]{\textbf{AMS subject classification. } #1}

\crefformat{equation}{(#2#1#3)}

\newcommand\nnfootnote[1]{%
  \begin{NoHyper}
  \renewcommand\thefootnote{}\footnote{#1}%
  \addtocounter{footnote}{-1}%
  \end{NoHyper}
}

\title{\normalfont{A Low-rank Approach for Nonlinear Parameter-dependent Fluid-structure Interaction Problems}} 

\author{ Peter Benner\textsuperscript{$\dag$}, Thomas Richter\textsuperscript{$\ddag$}, Roman Weinhandl\textsuperscript{$\S \star$}}

\date{} 

\begin{document}


\renewcommand{\sectionmark}[1]{\markright{A Low-rank Approach for FSI Problems}} 





\maketitle 







\begin{abstract}
Parameter-dependent discretizations of linear fluid-structure interaction problems can be approached with low-rank methods. When discretizing with respect to a set of parameters, the resulting equations can be translated to a matrix equation since all operators involved are linear. If nonlinear FSI problems are considered, a direct translation to a matrix equation is not possible. We present a method that splits the parameter set into disjoint subsets and, on each subset, computes an approximation of the problem related to the upper median parameter by means of the Newton iteration. This approximation is then used as initial guess for one Newton step on a subset of problems.
\end{abstract}



\nnfootnote{\textsuperscript{$\dag$} \textit{Max Planck Institute for Dynamics of Complex Technical Systems, Magdeburg and Otto von Guericke University Magdeburg, Germany} }
\nnfootnote{\textsuperscript{$\ddag$} \textit{Otto von Guericke University Magdeburg, Germany}}

\nnfootnote{\textsuperscript{$\S$} \textit{Otto von Guericke University Magdeburg and Max Planck Institute for Dynamics of Complex Technical Systems, Magdeburg, Germany (Email: \\weinhandl@mpi-magdeburg.mpg.de)}} 
\nnfootnote{\textsuperscript{$\star$}Corresponding author}

%




 \vspace{.4cm}
\par
\amssubj{15A69, 49M15, 65M22}

\section{Introduction}
Fluid-structure interaction (FSI) problems depend on parameters such as the solid shear modulus, the fluid density and the fluid viscosity. Parameter-\linebreak[4]dependent FSI discretizations allow to observe the reaction of an FSI model to a change of such parameters. A parameter-dependent discretization of a linear FSI problem yields many linear systems to be approximated. These equations can be translated to one single matrix equation. The solution, a matrix, can be approximated by a low-rank method as discussed in \cite{WeiBR19}. But as soon as nonlinear FSI problems are considered, such a translation is not possible anymore.
\par The proposed method extends the low-rank framework of \cite{WeiBR19} to nonlinear problems. It splits the parameter set into disjoint subsets. On each of these subsets, the Newton approximation for the problem related to the upper median parameter is computed to approximate the Jacobian matrix for all problems related to the subset. This allows to formulate a Newton step as a matrix equation. The Newton update, a matrix, can be approximated by a low-rank method and the global approximation to the parameter-dependent nonlinear FSI problem is achieved by stacking the approximations on the disjoint subsets column-wise.
\section{The Nonlinear Problem}
Let $d \in \{2,3\},$  $\Omega, F, S$ be open subsets of $\mathbb{R}^d$ with $\bar{F}\cup\bar{S}=\bar{\Omega}$, $F \cap S=\emptyset$. We use the stationary Navier-Stokes equations \cite[Section 2.4.5.3]{Ric17} to model the fluid part in $F$ and the stationary Navier-Lam\'e equations \cite[Problem 2.23]{Ric17} for the solid part in $S$. The interface is $\Gamma_{\text{int}}=\partial F \cap \partial S$, the boundary part where Neumann outflow conditions hold $\Gamma_f^\text{out}\subset \partial F \setminus \partial S$ and the boundary part where Dirichlet conditions hold $\Gamma_f^D=\partial F\setminus (\Gamma_f^{\text{out}} \cup \Gamma_{\text{int}})$. The weak formulation of the coupled nonlinear FSI problem with a vanishing right hand side $f$ reads
\begin{equation}\label{equation_weak_nonlinear1}
\begin{aligned}
\langle \nabla \cdot v, \xi \rangle_F&=0 \text{,}\\
\mu_s \langle \nabla u + \nabla u^T, \nabla \varphi \rangle_S+\lambda_s \langle \operatorname{tr}(\nabla u)I, \nabla \varphi \rangle_S\\
+\rho_s \langle (v \cdot \nabla) v, \varphi \rangle_F +\nu_f \rho_f \langle \nabla v + \nabla v^T, \nabla \varphi \rangle_F - \langle p, \nabla \cdot \varphi \rangle_F &=0 \quad \text{and}\\
\langle \nabla u, \nabla \psi \rangle_F &=0 \text{.}
\end{aligned}
\end{equation}
With $v_{\text{in}} \in H^1(\Omega)^d$, an extension of the Dirichlet data on $\Gamma_f^D$, the trial function $v \in v_{\text{in}}+H_0^1(\Omega, \Gamma_f^D \cup \Gamma_{\text{int}})^d$ is the velocity, $u \in H_0^1(\Omega)^d$ the deformation and $p \in L^2(F)$ the pressure. The test functions are $\xi \in L^2(F)$ (divergence equation), $\varphi \in H_0^1(\Omega, \partial \Omega \setminus \Gamma_f^{\text{out}})^d$ (momentum equation) and $\psi \in H_0^1(F)^d$ (deformation equation). The $\mathcal{L}^2$ scalar product on $F$ and $S$ is denoted by $\langle \cdot, \cdot \rangle_F$ and $\langle \cdot, \cdot \rangle_S$, respectively. The parameters involved are the kinematic fluid viscosity $\nu_f \in \mathbb{R}$, the fluid density $\rho_f \in \mathbb{R}$, the solid shear modulus $\mu_s \in \mathbb{R}$ and the first Lam\'e parameter $\lambda_s \in \mathbb{R}$.
\section{Discretization and Linearization}
Assume we are interested in discretizing the nonlinear FSI problem described in (\ref{equation_weak_nonlinear1})
parameter-dependently with respect to $m_1 \in \mathbb{N}$ shear moduli given by the set
\begin{align*}
S_\mu :=\{\mu_s^{i_1}\}_{i_1 \in \{1,...,m_1\}} \subset \mathbb{R}^+\text{,} \quad \text{with}\quad \mu_s^1 < ... < \mu_s^{m_1}\text{.}
\end{align*}
Consider a finite element discretization on $\Omega_h$, a matching mesh of the domain $\Omega$, with a total number of $N \in \mathbb{N}$ degrees of freedom. Let $A_0 \in \mathbb{R}^{N \times N}$ be a discretization matrix of all \textit{linear} operators involved in (\ref{equation_weak_nonlinear1}) with fixed parameters $\nu_f$, $\rho_f$, $\mu_s$ and $\lambda_s$. Let $A_1\in \mathbb{R}^{N \times N}$ be the discretization matrix of the operator
\begin{align}\label{operator_mu_lin1}
\langle \nabla u + \nabla u^T, \nabla \varphi \rangle_S.
\end{align}
The nonlinear part in (\ref{equation_weak_nonlinear1}), the convection term, requires a linearization technique.
\subsection{Linearization with Newton Iteration}
For a linearization by means of the Newton iteration, we need the Jacobian matrix of the operator $\langle (v \cdot \nabla)v,\varphi \rangle_F$. In our finite element space, every unknown $x_h=(p_h,v_h,u_h)^T \in \mathbb{R}^N$ consists of a pressure $p_h$, a velocity $ v_h$ and a deformation $u_h$. The discrete test space also has dimension $N$ and every unknown there can be written as $(\xi_h, \phi_h,\psi_h)^T \in \mathbb{R}^N$. The Jacobian matrix of $\langle (v \cdot \nabla)v,\varphi_h \rangle_F$ in our finite element space, evaluated at $x_h$, is
\small
\begin{align*}
J_{\langle (v \cdot \nabla)v,\varphi_h \rangle_F}(x_h)&=\left(  \begin{array}{ccc}
0 & 0 & 0 \\
0 & \frac{\partial \langle (v\cdot \nabla)v, \varphi_h \rangle_F}{\partial v}_{\big|_{v=v_h}}  &0 \\
0 & 0 & 0
\end{array}  \right) =:A_{\text{conv}}(x_h) \in \mathbb{R}^{N \times N} \text{.}
\end{align*}
\normalsize
\par Let $b_D \in \mathbb{R}^N$ be the right hand side vector that depends on the desired Dirichlet boundary conditions of the nonlinear FSI problem. Consider the FSI problem related to a fixed shear modulus $\mu_s^{i_1} \in S_\mu$ for some $i_1 \in \{1,...,m_1\}$ first.
\par If we start with an initial guess $x_0^{i_1}\in \mathbb{R}^N$, for instance $x_0^{i_1}=b_D$, 
at Newton step $j \in \mathbb{N}$, the equation
\begin{align}\label{equation_newton1}
\big(A_0+(\mu_s^{i_1}-\mu_s)A_1+\rho_f A_{\text{conv}}(x_{j-1}^{i_1})\big)s=-g(x_{j-1}^{i_1}, \mu_s^{i_1})
\end{align}
is to be solved for $s \in \mathbb{R}^N$. The approximation at linearization step $j$ then is
\begin{align*}
x_j^{i_1}=x_{j-1}^{i_1}+s \text{.}
\end{align*}
$g(x_{j-1}^{i_1},\mu_s^{i_1})$ evaluates all operators in (\ref{equation_weak_nonlinear1}) at the pressure, velocity and deformation of the approximation of the previous linearization step $x_{j-1}^{i_1}$ and the shear modulus $\mu_s^{i_1}$.
\section{Newton Iteration and Low-rank Methods}
In order to approximate a set of problems at one time, (\ref{equation_newton1}) has to be translated to a matrix equation. For this, first of all, we split the parameter set $S_\mu$ into disjoint subsets.
\par
If we perform a Newton step for a set of problems at one time, the same Jacobian matrix is used for the whole set. Therefore, the solutions to these different problems should not differ too much from each other. The method suggested in this paper splits the given parameter set into $K \in \mathbb{N}$ disjoint subsets, each of them containing adjacent parameters.
\begin{align*}
S_\mu= \bigcupdot \limits_{k=1}^K \mathcal{I}_k \text{.}
\end{align*}
\par By $\tilde{m}_k$, we denote the index of the upper median parameter of the set $\mathcal{I}_k$. After the parameter set is split into the subsets $\{\mathcal{I}_k\}_{k \in \{1,...,K\}}$,  we compute the Newton approximation $x_{\epsilon_N}^{\tilde{m}_k}$ of the problem related to the upper median parameter $\mu_s^{\tilde{m}_k}$ for all $k \in \{1,...,K\}$ up to some given accuracy $\epsilon_N>0$. $x_{\epsilon_N}^{\tilde{m}_k}$ is then used as initial guess for one Newton step.
\subsection{The Matrix Equation}  
 With $D_{1,k}:=\operatorname{diag}(\mathcal{I}_k)-\mu_s I^{|\mathcal{I}_k| \times |\mathcal{I}_k|}$ and $v_{\mathcal{I}_k}:=(\mu_s^{i_1})_{i_1\in \mathcal{I}_k} \in \mathbb{R}^{|\mathcal{I}_k|}$, the matrix equation that is to be solved for $S_k \in \mathbb{R}^{N \times |\mathcal{I}_k|}$ on every subset $\mathcal{I}_k$ is
\small
\begin{equation}\label{equation_global_newton1}
\begin{aligned}
A_0S_k+A_1S_kD_{1,k}+\rho_f A_{\text{conv}}(x_{\epsilon_N}^{\tilde{m}_k})S_k=&\underbrace{-g(x_{\epsilon_N}^{\tilde{m}_k}, 0) \otimes (1,...,1)-A_1 x_{\epsilon_N}^{\tilde{m}_k} \otimes  v_{\mathcal{I}_k}^T}_{=:B_k} \text{.}
\end{aligned}
\end{equation}
\normalsize
$I^{|\mathcal{I}_k| \times |\mathcal{I}_k|}$ denotes the identity matrix of size $|\mathcal{I}_k| \times |\mathcal{I}_k|$. In (\ref{equation_global_newton1}), the initial guess for the Newton step is
\begin{align*}
X_{\text{initial}}^k:=x_{\epsilon_N}^{\tilde{m}_k} \otimes (1,...,1)\text{.}
\end{align*}
The approximation at the next linearization step is 
\begin{align}\label{equation_newton_app_linstep2}
X^k:=X_{\text{initial}}^k+S_k\text{.}
\end{align}
\par The global approximation for the whole parameter-dependent problem then is
\begin{align*}
\tilde{X}:=[X^1|\cdots|X^K]\text{.}
\end{align*}
\begin{remark} The initial guess for the Newton step (\ref{equation_global_newton1}), $X_{\text{initial}}^k$, has rank $1$ and the operator (\ref{operator_mu_lin1}) is linear. This is why the rank of the right hand side matrix $B_k$ in (\ref{equation_global_newton1}) is not bigger than $2$.
\end{remark}
\begin{remark}
If multiple Newton steps like (\ref{equation_global_newton1}) were performed, two main difficulties would come up. At step $2$, the approximation of the previous linearization step would be given by $X^k$ from (\ref{equation_newton_app_linstep2}).
\par \textbf{The right hand side:} $X^k$ is not a matrix of low rank and $g(\cdot,\cdot)$ would have to be evaluated for all columns of $X^k$ separately in a second Newton step. Thus, the right hand side matrix $B_k$ would not have low-rank structure either.
\par \textbf{The Jacobian matrix: }Since all columns of the initial guess $X_{\text{initial}}^k$ coincide, the Jacobian matrix in (\ref{equation_global_newton1}) is correct for all equations related to the parameter set $\mathcal{I}_k$. But the columns of $X^k$ differ from each other. A second Newton step would then become what is, in the literature, often called an inexact Newton step \cite[Remark 5.7]{Ric17}.
\end{remark}
\subsection{Low-rank Methods}
Let $k \in \{1,...,K\}$,
\begin{align*}
\tilde{A}_{\text{conv}}:=A_{\text{conv}}(x_{\epsilon_N}^{\tilde{m}_k}) \quad \text{and}\quad b_g:=g(x_{\epsilon_N}^{\tilde{m}_k},0)\text{.}
\end{align*}
Consider only the column related to the parameter index $i_1\in \mathcal{I}_k$ in (\ref{equation_global_newton1}):
\begin{align*}
\underbrace{\big(A_0+(\mu_s^{i_1}-\mu_s)A_1+\rho_f \tilde{A}_{\text{conv}}\big)}_{=:A(\mu_s^{i_1})}s^{i_1}=\underbrace{-b_g-\mu_s^{i_1}A_1 x_{\epsilon_N}^{\tilde{m}_k}}_{=:b(\mu_s^{i_1})} \text{,} \qquad \text{with} \qquad s^{i_1} \in \mathbb{R}^N\text{.}
\end{align*}
Assume that $x_{\epsilon_N}^{\tilde{m}_k}$ is fixed and $A(\mu_s^{i_1})$ is invertible for all $\mu_s^{i_1} \in \mathcal{I}_k$. $A(\mu_s^{i_1})$ and $b(\mu_s^{i_1})$ depend linearly on $\mu_s^{i_1}$. $A(\cdot)$ and $b(\cdot)$ are analytic matrix- and vector-valued functions, respectively. Due to \cite[Theorem 2.4]{KreT11}, the singular value decay of the matrix $S_k$ in (\ref{equation_global_newton1}) is exponential. Algorithm \ref{algorithm1_lowrank} exploits this fact and approximates $S_k$ in (\ref{equation_global_newton1}) by a low-rank matrix.
\section{Numerical Results}
A $3d$ jetty flow in a channel with the geometric configuration
\begin{align*}
\Omega:=(0,12) \times (0,4) \times (0,4) \text{, } S:=(2,3) \times (0,2) \times (0,4) \quad  \text{and} \quad  F:=\Omega \setminus \bar{S} 
\end{align*}
is considered. The left Dirichlet inflow is given by the velocity profile
\begin{align*}
v=\left(\begin{array}{c}v_1\\v_2\\v_3
\end{array} \right) =\left( \begin{array}{c}
 \frac{1}{10}y(4-y)z(4-z)\\0\\0 \end{array} \right) \in \mathbb{R}^3 \quad \text{at} \quad x=0\text{.}
\end{align*}
At $x=12$, the do nothing boundary condition holds. At $z=0$, deformation and velocity in  normal direction is prohibited. Everywhere else on $\partial \Omega$, the velocity and the deformation vanish. For the Navier-Stokes equations, stabilized Stokes elements \cite[Lemma 4.47]{Ric17} are used.
\begin{algorithm}[h]
\caption{Low-rank Method for One-parameter Nonlinear FSI}
\begin{algorithmic}\label{algorithm1_lowrank}
\REQUIRE{Accuracy $\epsilon_N>0$ for Newton method, ranks $R_k \in \mathbb{N}$ for $k \in \{1,...,K\}$}
\ENSURE{The rank-$\sum \limits_{k=1}^K R_k$ approximation $\hat{X}$ of the parameter-dependent FSI discretization}
\STATE Split the parameter set $S_\mu$ into the subsets $\bigcupdot \limits_{k=1}^K \mathcal{I}_k$.
\FOR{$k=1,...,K$}
\STATE Compute the Newton approximation of the upper median parameter problem related to a shear modulus of $\mu_s^{\tilde{m}_k}$ with accuracy $\epsilon_N$ $\curvearrowright$ $x_{\epsilon_N}^{\tilde{m}_k}$.
\STATE Use $x_{\epsilon_N}^{\tilde{m}_k} \otimes (1,...,1) \in \mathbb{R}^{N \times |\mathcal{I}_k|}$ as initial guess for one Newton step on $\mathcal{I}_k$. Find a rank-$(R_k-1)$ approximation $\hat{S}_k$ that approximates $S_k \in \mathbb{R}^{N \times |\mathcal{I}_k|}$ from (\ref{equation_global_newton1}) by a low-rank method from \cite{WeiBR19}. 
\STATE Build the sum
\begin{align*}
\hat{X}_k=x_{\epsilon_N}^{\tilde{m}_k} \otimes (1,...,1)+\hat{S}_k\text{.}
\end{align*}
\ENDFOR
\STATE $\hat{X}:=[\hat{X}_1 | ... | \hat{X}_K]$
\end{algorithmic}
\end{algorithm}
\subsection{Parameters}
The nonlinear FSI problem (\ref{equation_weak_nonlinear1}) is discretized with $Q_1$ elements (compare \cite[(4.13)]{Ric17}) with respect to 
\begin{align*}
1500 \text{ shear moduli } \mu_s^{i_1} \in S_\mu \subset  [20000,60000] \text{.}
\end{align*}
The fixed first Lam\'e parameter is $\lambda_s=200000$. With these parameters, solid configurations with Poisson ratios between $0.38462$ and $0.45455$ are covered. The fluid density is $\rho_f=12.5$ and the kinematic fluid viscosity is $\nu_f=0.04$.
\subsection{Comparison ChebyshevT with Standard Newton}A server operating CentOS 7 with 2 AMD EPYC 7501 and 512GB RAM, MATLAB\textsuperscript{\textregistered} 2017b in combination with the htucker MATLAB toolbox \cite{software_htucker1} and the finite element toolkit GASCOIGNE \cite{software_gascoigne1} was used to compare Algorithm \ref{algorithm1_lowrank} with $1500 $ Newton iterations applied consecutively. The parameter set $S_\mu$ was split into $K=15$ subsets.
\begin{figure}[ht]
\centering
\includegraphics[scale=.57]{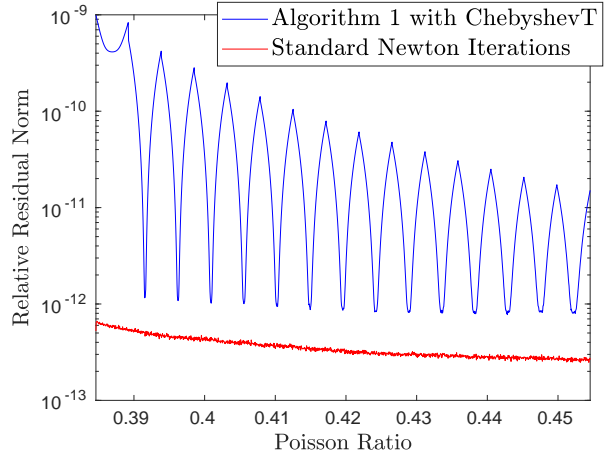}
\caption{Comparison of the approximations provided by Algorithm \ref{algorithm1_lowrank} and $1500$ standard Newton iterations applied consecutively.}
\label{figure_plot1}
\end{figure}
\subsubsection{Preconditioner and Eigenvalue Estimation}
After the Newton approximations $x_{\epsilon_N}^{\tilde{m}_k}$ are available for all $k \in \{1,...,15\}$, the LU decomposition of the mean-based preconditioner $P_T^k$ \cite[Section 3.2]{WeiBR19} of $A(\mu_s^{i_1})$ is computed separately on every subset $\mathcal{I}_k$. To estimate the parameters $d$ and $c$ for the ChebyshevT method \cite[Algorithm 3]{WeiBR19}, the eigenvalues of the matrices
\begin{align*}
(P_T^{k})^{-1} A\Big(\max \limits_{i_1 \in \mathcal{I}_k}(\mu_s^{i_1})\Big) \qquad \text{and} \qquad (P_T^k)^{-1} A\Big(\min \limits_{i_1 \in \mathcal{I}_k}(\mu_s^{i_1}) \Big)
\end{align*}
are taken into consideration. For all $K=15$ subsets, they were estimated to $d=1$ and $c=0.1$ for a small number of $N=945$ degrees of freedom within $28.3$ seconds (computation time for the Newton approximations included). Every approximation is related to a certain shear modulus. Therefore, all problems differ by the Poisson ratio of the solid. The y-axis in Figure \ref{figure_plot1}, on the other hand, corresponds to the relative residual norm \begin{align*}
\frac{\|g(x^{i_1},\mu_s^{i_1})\|_2}{\|g(b_D,\mu_s^{i_1})\|_2}
\end{align*}
of the approximation $x^{i_1}$ for $i_1 \in \{1,...,m_1\}$. Algorithm \ref{algorithm1_lowrank} was applied with $\epsilon_N=10^{-4}$ and $R_k=10$ $\forall k \in \{1,...,15\}$ to a problem with $N=255255$. Therefore, the global approximation rank is $R=150$. In comparison to this, standard Newton iterations were applied to the $1500$ separate problems consecutively where for every Newton iteration, the last approximation served as initial guess for the next Newton iteration.
%
\par The approximations obtained by the Standard Newton iterations within $238$ hours ($1507$ Newton steps) provided, as visualized in Figure \ref{figure_plot1}, relative residuals with norms smaller than $10^{-12}$ each. Algorithm \ref{algorithm1_lowrank} took $519$ minutes ($35$ Newton steps) to compute the low-rank approximation. In addition to the $28.3$ seconds for the eigenvalue estimation, the $20$ Newton steps to compute $x_{\epsilon_N}^{\tilde{m}_k}$ for all $k \in \{1,...,15\}$ took, in total, $195.65$ minutes and the $15$ Newton steps for the matrix equations (\ref{equation_global_newton1}) took, in total, $323.3$ minutes.
\section{Conclusions}
Low-rank methods can be used to compute approximations to parameter-\linebreak[4]dependent nonlinear FSI discretizations, in particular, if each of the subsets, the parameter set is split into, does contain problems that do not differ too much from each other. The Newton step on the subset uses the same Jacobian matrix and the same initial guess for the whole subset. It has to provide acceptable convergence within one single step not only for the upper median problem.
\par Whether the results can be improved by choosing the subsets $\mathcal{I}_k$ or the approximation ranks on these subsets adaptively, is still open. Moreover, how these low-rank methods can be applied to fully nonlinear FSI problems that use, in addition to the Navier-Stokes equations on the fluid, for instance, the St. Venant Kirchhoff model equations \cite[Definition 2.18]{Ric17} on the solid is an open problem. Then, the right hand side in (\ref{equation_global_newton1}) would have to be approximated.

\section*{Acknowledgements}This work was supported by the Deutsche Forschungsgemeinschaft (DFG, German Research Foundation) - 314838170, GRK 2297 MathCoRe.



\bibliographystyle{siam}



\end{document}